\documentclass[12pt,a4paper]{article}

\usepackage{float}
\usepackage[left=2cm,right=2cm,top=2cm,bottom=2cm]{geometry}
\usepackage[dvips]{graphics}
\usepackage[dvips]{epsfig}
\usepackage{color}
\usepackage{hyperref}
\usepackage{tikz,authblk}

\newtheorem{lemma}{Lemma}
\newtheorem{theorem}[lemma]{Theorem}
\newtheorem{corollary}[lemma]{Corollary}
\newtheorem{proposition}[lemma]{Proposition}
\newtheorem{definition}{Definition}
\newtheorem{remark}{Remark}

\newtheorem{example}{Example}

\newcommand{\dimo}{\noindent \emph{Proof. }}
\newcommand{\qed}{\\ \rightline{$\Box$}\\}
\newcommand{\e}{\varepsilon}
\newcommand{\G}{\Gamma}
\newcommand{\g}{\gamma}

\usepackage{latexsym}
\usepackage{amsfonts}
\usepackage{amssymb}
\usepackage{amsmath}
\usepackage[english]{babel}

\begin{document}

\title{Estimating trisection genus via gem theory}

 \renewcommand{\Authfont}{\scshape\small}
 \renewcommand{\Affilfont}{\itshape\small}
 \renewcommand{\Authand}{ and }

\author[1] {Maria Rita Casali}
\author[2] {Paola Cristofori}

\affil[1] {Department of Physics, Mathematics and Computer Science, University of Modena and Reggio Emilia, 
Via Campi 213 B, I-41125 Modena (Italy), casali@unimore.it}

\affil[2] {Department of Physics, Mathematics and Computer Science, University of Modena and Reggio Emilia, Via Campi 213 B, I-41125 Modena (Italy), paola.cristofori@unimore.it}

\maketitle
 
 \begin{abstract}
{\it Gems} are a particular type of edge-colored graphs, dual to colored triangulations, which represent compact PL-manifolds of arbitrary dimension, both in the closed and boundary case. 
In the present paper, gem theory is used to approach trisections of PL 4-manifolds, so as to prove that: 
\begin{itemize}
  \item the graph-defined invariant {\it regular genus} is an upper bound for the trisection genus of each closed $4$-manifold; 
  \item a trisection diagram can be directly obtained from any gem of a closed $4$-manifold. 
\end{itemize}
Moreover, suitable extensions of the above results are presented for compact 4-manifolds with connected boundary. 
\end{abstract}

\maketitle

\medskip



\section{Introduction}

The notion of trisection of a smooth closed 4-manifold (originally introduced by Gay and Kirby in \cite{Gay-Kirby}) is well-known: each such manifold can be decomposed into three 4-dimensional handlebodies, with disjoint interiors and mutually intersecting in 3-dimensional handlebodies (with the same genus), so that the intersection of all three ``pieces" is a closed surface. The obvious analogy with Heegaard splittings in dimension three suggests also to define an associated invariant in dimension four, i.e. the {\it trisection genus}, as the minimum genus of the 3-dimensional handlebodies involved in all trisections of the given 4-manifold.

On the other hand, {\it gem theory} enables to combinatorially represent PL-manifolds of arbitrary dimension via edge-colored graphs, both in the closed and boundary case; within this theory, the invariant {\it regular genus} has been defined and extensively studied, and classification results have been obtained according to it. 

The present paper follows some previous works, where the relationship between the two theories has been established, allowing to obtain trisections of closed 4-manifolds from gems, and also to suitably extend the notion of trisection to 4-manifolds with non-empty connected boundary (by making use of the notion of {\it gem-induced trisection}, which involves an Heegaard splitting of the boundary). 
However, the existing constructions 
either are limited to  gems satisfying certain conditions (\cite{Spreer-Tillmann(Exp)}, \cite{Casali-Cristofori gem-induced} and \cite{Casali-Cristofori trisection bis}), or give rise to trisections with a very high genus (\cite{Martini-Toriumi}). 

In this paper, we show how to obtain, from any 
gem of a compact 4-manifold with empty (resp. connected) boundary, a trisection (resp. a gem-induced trisection), while keeping under control the genus of the splitting surface. 

As a consequence, the following estimations of the trisection genus are proved: 

\begin{theorem}\label{intro:theorem}
Let $M$ be a compact $4$-manifold whose boundary is either empty or a connected sum of sphere bundles over $\mathbb S^1$; then: 
$$g_T(\bar M) \le \mathcal G(M)$$
\noindent where $\mathcal G(M)$ denotes the regular genus of $M$, while 
$g_T(\bar M)$ is the trisection genus of the closed 4-manifold $\bar M$ (uniquely) associated to $M$.  

\medskip 

\noindent In particular, if $M$ is closed, then: $$g_T(M) \le \mathcal G(M).$$
\end{theorem}

 Hence, the trisection genus of a closed 4-manifold turns out to admit as an upper bound, not only its regular genus, 
but also that of the compact 4-manifold consisting of handles up to index two in a handle-decomposition of the closed 4-manifold itself. 
The estimation is established by constructing explicitly a trisection of the closed manifold via gems.
 
At the same time, this construction enables one to obtain a trisection diagram from any gem of the closed $4$-manifold, by making use of its bicolored cycles: see Proposition \ref{trisection_diagrams_gem-induced(chiuse)}. 
\bigskip

It should be noted that the whole work is carried out for compact 4-manifolds with connected boundary, by means of the concept of gem-induced trisection: 
closed manifolds are considered as a particular case through their identification with manifolds with spherical boundary.




More precisely, if $M$ is a compact 4-manifold with connected boundary, we get an upper and lower estimation of the so called {\it G-trisection genus} of $M$ (the analogue of trisection genus for gem-induced trisections: see Definition \ref{def_GT-genus}) in terms of the regular genus of $M$, of the Heegaard genus of its boundary $\partial M$ and of the rank of the fundamental group of the associated singular manifold, obtained by capping off $\partial M$ by a cone: see Corollary \ref{trisection_vs_regular-genus}.  

Moreover, if either $\partial M$ is a connected sum of sphere bundles over $\mathbb S^1$ or $M$ is simply-connected, starting from each gem of $M$ we describe how to obtain a so called {\it G-trisection diagram} of $M$ (a suitable generalization of trisection diagrams: see Definition \ref{def. G-trisection-diagram}), consisting of bicolored cycles of the gem itself: see Proposition \ref{trisection_diagrams_gem-induced(chiuse)} and Proposition \ref{trisection_diagrams_gem-induced(bordo)}. 

\bigskip 
The paper is structured as follows: after a brief review on gem theory (Section \ref{ss:gem-theory}), existing results on trisections arising from gems are recalled in Section \ref{s:gem-induced-trisections}, together with the extension of the notion of {\it stabilization} in order to ensure the existence of gem-induced trisections also from gems of non-orientable manifolds and/or of compact manifolds with non-empty connected boundary; then, Section \ref{ss:proving-estimations} is devoted to prove the estimation results concerning trisection genus and G-trisection genus, while trisection diagrams and G-trisection diagrams arising from gems are the subject of Section \ref{ss.from_gems_to_trisection_diagrams}.     

\bigskip 

\section{Basic notions of gem theory}
\label{ss:gem-theory} 

The present section is devoted  to a brief  review of some basic notions of the so called  {\it gem theory}, which is a representation tool for piecewise linear (PL) compact manifolds, without assumptions about dimension, connectedness, orientability or boundary properties  (see the ``classical'' survey paper \cite{Ferri-Gagliardi-Grasselli}, or the more recent one \cite{Casali-Cristofori-Gagliardi Complutense 2015}, concerning the $4$-dimensional case). 

From now on, unless otherwise stated, all spaces and maps will be considered in the PL category, and all manifolds will be assumed to be compact and connected.  

\begin{definition} \label{$n+1$-colored graph} {\rm An $(n+1)${\emph{-colored graph}}  ($n \ge 2$) is a pair $(\G,\g)$, where $\G=(V(\G), E(\G))$ is a multigraph (i.e. it can contain multiple edges, but no loops) 
which is regular of degree  $n+1$, and $\g$ is an {\it edge-coloration}, that is a map  $\g: E(\G) \rightarrow \Delta_n=\{0,\ldots, n\}$ which is injective on 
adjacent edges.}
\end{definition}

In the following, for the sake of conciseness, 
when the coloration is clearly understood, we will denote the colored graph simply by $\G$. 

\smallskip

For every  $\{c_1, \dots, c_h\} \subseteq\Delta_n,$ let $\G_{\{c_1, \dots, c_h\}}$  be the subgraph obtained from $\G$  by deleting all 
 edges that are not colored by the elements of $\{c_1, \dots, c_h\}$. 
In this setting, the complementary set of $\{c\}$ 
in $\Delta_n$ will be denoted by $\hat c$. 
The connected components of $\G_{\{c_1, \dots, c_h\}}$ are called {\it $\{c_1, \dots, c_h\}$-residues} or {\it $h$-residues} of $\G$; their number is denoted by $g_{\{c_1, \dots, c_h\}}$ (or, for short, by $g_{c_1,c_2}$, $g_{c_1,c_2,c_3}$ and $g_{\hat c}$ if $h=2,$ $h=3$ and $h = n$ respectively). 

 \medskip 

\noindent An $n$-dimensional pseudocomplex $K(\G)$ can be associated to any $(n+1)$-colored graph $\G$ as follows: 
\begin{itemize}
\item for each vertex of $\G$, take an $n$-simplex and label its vertices by the elements of $\Delta_n$;
\item whenever two vertices of $\G$ are $c$-adjacent ($c\in\Delta_n$), glue the corresponding $n$-simplices  along their $(n-1)$-dimensional faces opposite to the $c$-labeled vertices, so as to identify equally labeled vertices.
\end{itemize}

\smallskip
In general $\vert K(\G)\vert$ turns out to be an {\it n-pseudomanifold} and $\G$ is said to {\it represent} it.  
 
\medskip

Note that, by construction, $K(\G)$ is endowed with a vertex-labeling by $\Delta_n$ that is injective on each 
simplex. Moreover, $\G$ can be identified with 
the 1-skeleton of the dual complex of $K(\G)$.
The duality establishes a bijection between the $\{c_1, \dots, c_h\}$-residues of  $\G$  
and the $(n-h)$-simplices of $K(\G)$ whose vertices are labeled by   $\Delta_n \setminus \{c_1, \dots, c_h\}$. 

\noindent 

In particular, if $\G$ is an $(n+1)$-colored graph, each connected component of $\G_{\hat c}$ ($c\in\Delta_n$) is an $n$-colored graph representing the link of a $c$-labeled vertex of $K(\G)$ in the first barycentric subdivision of $K(\Gamma).$

Therefore: 
\begin{itemize}
\item $\vert K(\G)\vert$ is a {\it closed $n$-manifold} iff, for each color $c\in\Delta_n$, all $\hat c$-residues of $\G$ represent the $(n-1)$-sphere; 
\item $\vert K(\G)\vert$ is a {\it singular\footnote{A {\it singular (PL) $n$-manifold} is a closed connected $n$-dimensional polyhedron admitting a simplicial triangulation where the links of vertices are closed connected $(n-1)$-manifolds, while
the links of the $h$-simplices, with $h > 0$, are $(n-h-1)$-spheres. Vertices whose links are not PL $(n-1)$-spheres are called {\it singular}.

The notion extends also to polyhedra associated to colored graphs, provided that links 
are considered in the first barycentric subdivision: 
in this case, the condition regarding $h$-simplices, $h >0$, directly follows from that about 
vertices.}   
$n$-manifold} iff, for each color $c\in\Delta_n$, all $\hat c$-residues of $\G$ represent closed connected $(n-1)$-manifolds.
\end{itemize}

\begin{remark} \label{correspondence-sing-boundary} {\em Given a compact $n$-manifold $M$, a singular $n$-manifold $\widehat M$ can be constructed by capping off each component of $\partial M$ by a cone over it.
Conversely, if $N$ is a singular $n$-manifold, then a compact $n$-manifold is easily obtained by deleting small open neighbourhoods of its singular vertices.
Note that, as a consequence, compact $n$-manifolds with no spherical boundary components and singular $n$-manifolds correspond bijectively  to each other in a well-defined way.   
For this reason and since compact manifolds with spherical boundary components may be identified with the associated closed ones, throughout this work, we will restrict our attention to compact manifolds without spherical boundary components. Obviously, in this context, closed $n$-manifolds are characterized by $M= \widehat M.$ }
\end{remark}

The bijection described in Remark \ref{correspondence-sing-boundary} allows to say that an $(n+1)$-colored graph $\G$ {\it represents}
a compact $n$-manifold $M$ with no spherical boundary components (or, equivalently, is 
 a {\it gem} of $M$, where gem means {\it Graph Encoding Manifold})  if and only if  it represents the associated singular manifold $\widehat M$.

If $\G$ is a gem of a compact $n$-manifold $M$,  an $n$-residue of $\G$ will be called {\it singular} if it does not represent $\mathbb S^{n-1}$; similarly, a color $c$ will be called {\it singular} if at least one of the $\hat c$-residues of $\G$ is singular. Note that each singular $n$-residue turns out to represent a connected component of $\partial M.$

\medskip
The following theorem extends to the boundary case a well-known result 
(originally due to Pezzana, as explained in the cited surveys), 
which is the basis of the combinatorial representation theory for closed manifolds of arbitrary dimension via colored graphs. 

\begin{theorem}{\em (\cite{Casali-Cristofori-Grasselli})}\ \label{Theorem_gem}  
Any compact orientable (resp. non orientable) $n$-manifold with no spherical boundary components admits a bipartite (resp. non-bipartite) $(n+1)$-colored graph representing it.
\par \noindent In particular, if  $M$ has empty or connected boundary: 
\begin{itemize}
\item $\G$ may be assumed to have color $n$ as its unique possible singular color, and exactly one $\hat n$-residue 
(in this case we will say that $\G$ belongs to the class $G^{(n)}_s$);\footnote{Equivalently, the colored triangulation $K(\G)$ of $\widehat M$  has exactly one vertex labelled $n$,  which is the only possible singular vertex of $K(\G)$. 
Note that the vertex of $K(\G)$ labelled $n$ (or, equivalently, the $\hat n$-residue of $\G$) is singular if and only if $\partial M$ is non-empty and non-spherical. 
} 
\item $\G$ may be assumed to have exactly one $\hat c$-residue, $\forall c \in \Delta_n$ (in this case $\G$ is  called  
a \emph{crystallization} of $M$).\footnote{Equivalently,  $K(\G)$ 
 has exactly $n+1$ vertices.}  
\end{itemize}
\end{theorem}

\bigskip
The notion of regular genus - which is a PL invariant of great importance within gem theory - is based on the existence of a particular type of embedding of colored graphs into surfaces. 

\begin{proposition}{\em (\cite{Ferri-Gagliardi-Grasselli})}\label{reg_emb}
Let $\G$ be a connected bipartite  (resp. non-bipartite)
$(n+1)$-colored graph of order $2p$. Then for each cyclic permutation $\varepsilon = (\varepsilon_0,\ldots,\varepsilon_n)$ of $\Delta_n$, up to inverse, there exists a cellular embedding, called \emph{regular}, of $\G$  
into an orientable (resp. non-orientable)
closed surface $F_{\varepsilon}(\G)$ whose regions are bounded by the images of the $\{\varepsilon_j,\varepsilon_{j+1}\}$-colored cycles, for each $j \in \mathbb Z_{n+1}$.
Moreover, the genus (resp. half the genus)
$\rho_{\varepsilon} (\G)$  of $F_{\varepsilon}(\G)$ satisfies

\begin{equation*}
2 - 2\rho_\varepsilon(\G)= \sum_{j\in \mathbb{Z}_{n+1}} g_{\varepsilon_j, \varepsilon_{j+1}} + (1-n)p.
\end{equation*}
\end{proposition}

\begin{definition} {\rm The \emph{regular genus} of  an $(n+1)$-colored graph $\G$ is defined as
$$\rho(\G) = min\{\rho_\varepsilon(\G)\ \vert \ \varepsilon\ \text{cyclic permutation of \ } \Delta_n\}; $$
 the  {\it regular genus} of a compact $n$-manifold $M$ is defined as
$$\mathcal G (M) = min\{\rho(\G)\ \vert \ \G\ \text{gem of \ } M\}.$$}
\end{definition}

\medskip
 
\begin{remark} \label{rem_regular genus} {\em  The regular genus extends to higher dimension the classical genus of a surface and the Heegaard genus of a $3$-manifold. 
In addition to the characterization of spheres in any dimension, 
it yielded  many classification results, 
especially in dimension $4$ and $5$ (see 
\cite{Casali-Cristofori-Gagliardi Complutense 2015}, \cite{generalized-genus}  and their references). Moreover, the regular genus is strictly related to the {\it G-degree}, a PL invariant arising within the theory of {\it Colored tensor models} in theoretical physics (\cite{Casali-Cristofori-Dartois-Grasselli}, \cite{Casali-Cristofori-Grasselli}, \cite{Casali-Grasselli 2019}). }
\end{remark}





\bigskip 

\section{Trisections arising from gems} \label{s:gem-induced-trisections} 
The notion of trisection of a smooth closed orientable $4$-manifold has been introduced by Gay and Kirby in 2016 and more recently extended to the non-orientable case by Miller and Naylor. The following definition, which is commonly used at present, is a slight generalization of Gay and Kirby's. 
 
\begin{definition}\label{def: trisection}{\em  \ {\rm (\cite{Gay-Kirby}, \cite{Spreer-Tillmann(Exp)}, \cite{Miller-Naylor})}
A {\it trisection of genus $g$} of a smooth closed orientable (resp. non-orientable) $4$-manifold $M$ is a decomposition of $M$ into three 4-dimensional orientable (resp. non-orientable) handlebodies $H_0,\ H_1, \ H_2$,  with disjoint interiors, such that their  pairwise intersections are  $3$-dimensional orientable (resp. non-orientable) handlebodies  of genus $g$ and  $H_{0}\cap H_{1}\cap H_{2}$  is a closed connected orientable (resp. non-orientable) surface. \\
The {\it trisection genus} $g_T(M)$ of $M$ is defined as the minimum genus of all 
trisections of $M.$
 }
\end{definition} 

\begin{remark}\label{rem:spine}{\em A well-known result by Laudenbach and Poenaru for the orientable case (\cite{Laudenbach-Poenaru}) and its non-orientable counterpart by Miller and Naylor (\cite{Miller-Naylor}) allow to prove that 
the manifold is completely determined by the union of the three 3-dimensional handlebodies of any trisection  (see also \cite{Meier-Schirmer-Zupan}).
}
\end{remark}

The first link between trisections and gems was established by Spreer and Tillmann (\cite{Spreer-Tillmann(Exp)}) by exploiting Bell-Hass-Rubinstein-Tillmann's approach to trisections through tricolorings of triangulations (\cite{Bell-et-al}) and by applying it to colored triangulations whose 1-skeleton equals that of a 4-simplex.

\medskip

More generally, it is shown in  \cite{Casali-Cristofori gem-induced} and \cite{Casali-Cristofori trisection bis} that, starting from a gem $\G\in G_s^{(4)}$ of a compact $4$-manifold $M$ with empty or connected boundary and a cyclic permutation $\varepsilon=(\e_0,\e_1,\e_2,\e_3,\e_4=4)$ of $\Delta_4$, it is always possible to construct a decomposition $\mathcal  T(\Gamma, \varepsilon)=(H_0,H_1,H_2)$ of $M$ into three submanifolds with  pairwise disjoint interiors such that:
 \begin{itemize}
 \item [-]  $H_{1},H_{2}$ are $4$-dimensional handlebodies, while $H_{0}$ is a $4$-disk or a collar of $\partial M$ according to the boundary being empty or not;
\item [-] the intersections $H_{01}=H_{0}\cap H_{1}$ and  $H_{02}=H_{0}\cap H_{2}$ are $3$-dimensional handlebodies of genus $\rho_{\varepsilon}(\G_{\hat 4})$.
 \item [-] $\Sigma = H_{0}\cap H_{1}\cap H_{2}$  is a closed connected surface ({\it central surface} of  $\mathcal  T(\Gamma, \varepsilon)$).  Moreover, $(\Sigma, H_{01}, H_{02})$ is a (genus $\rho_{\varepsilon}(\G_{\hat 4})$) Heegaard splitting of $|K(\G_{\hat 4})|$.
 \end{itemize}
 
If $H_{12}=H_1\cap H_2$ is a 3-dimensional handlebody, too, then $\mathcal  T(\Gamma, \varepsilon)$ is called a {\it gem-induced trisection} of $M$
\footnote{We point out that, while in the closed case gem-induced trisections trivially satisfy Definition \ref{def: trisection}, in the case of non-empty boundary a gem-induced trisection is not a trisection in the sense of  \cite{Castro-Gay-Pinzon}: 
for instance, it determines a Heegaard splitting of the boundary  instead of an open-book decomposition
.}. Similarly to trisection theory, the common genus of all 3-dimensional handlebodies is called the {\it genus} of $\mathcal  T(\Gamma, \varepsilon)$.

Moreover, if  $\mathcal  T(\Gamma, \varepsilon)$ is a gem-induced trisection, then all the involved handlebodies, as well as the central surface $\Sigma$, are orientable or not according to the orientability of $M$. 

\medskip 

The construction of $\mathcal  T(\Gamma, \varepsilon)=(H_0,H_1,H_2)$ makes use of a tricoloring of the vertices of $K(\G)$ giving color red to the $\e_0$- and $\e_2$-labelled vertices, color green to the $\e_1$- and $\e_3$-labelled vertices, and color blue to the only vertex labelled by $\e_4=4.$\footnote{In the closed case, the role of color $4$ can be played by any color $c$ such that $\Gamma_{\hat c}$ is connected.} 

\begin{figure}   [h!]
     \centering
    \includegraphics[width=1\linewidth]{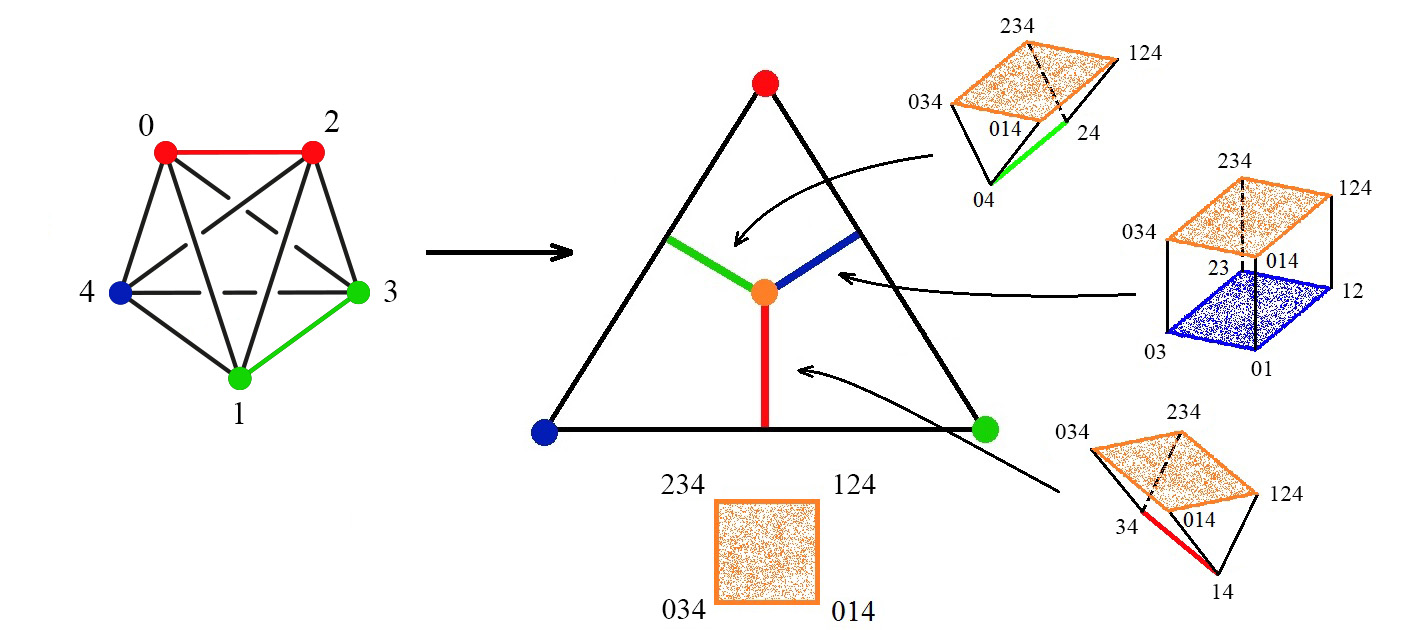}
    \caption{The intersections of the 3-dimensional pieces and the central surface of $\mathcal  T(\Gamma, \varepsilon)$ with any 4-simplex  (redrawn from \cite{Spreer-Tillmann(Exp)})}
     \label{fig:trisec_pieces}
 \end{figure}

The preimage of the cubical decomposition of the standard 2-simplex under the simplicial map induced by the tricoloring gives rise to a decomposition  of $|K(\G)|=\widehat M$ into a triple 
$(\widehat H_0,H_1,H_2)$, where the $\widehat H_0$ is a regular neighbourhood of the blue vertex, while $H_1$ and $H_2$ are regular neighbourhoods of the 1-dimensional subcomplexes of $K(\Gamma)$ generated by the red and green vertices respectively (see Figure \ref{fig:trisec_pieces}). By deleting from $\widehat H_0$ a suitable small neighbourhood of the blue vertex, the triple $(H_0,H_1,H_2)$ is finally obtained  (see \cite{Casali-Cristofori gem-induced}, \cite{Casali-Cristofori trisection bis} for details).

The 3-dimensional intersections $H_{01}, H_{02}$ and  $H_{12}$ meet each 4-simplex in  two prisms and one cube respectively, whose vertices are barycenters of 1- and 2-simplexes of the triangulation: see the right side of Figure \ref{fig:trisec_pieces}, where these barycenters are indicated by the labels of the spanning vertices of their corresponding simplexes (for simplicity $\e_i=i,\ \forall i\in\Delta_4$ is also assumed). 

As a consequence, $H_{01}$ (resp. $H_{02}$) is easily seen to be a 3-dimensional handlebody since each of its constituting prisms collapses to the green (resp. red) edge in the figure from its opposite face, which, as part of the surface $\Sigma$, is free; by the same reason, there is an obvious collapse of $H_{12}$ to the 2-dimensional complex $Q(\G,\varepsilon)$ formed by the squares that are bottom faces of its constituting cubes (i.e. by the squares depicted in blue in Figure \ref{fig:trisec_pieces}).

\begin{figure}[h!]
\centering
\scalebox{0.6}{\includegraphics{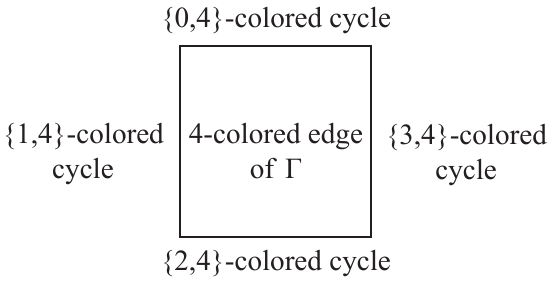}}
\caption{the square, corresponding to a $4$-colored edge of $\G$, constituting $Q(\G,\varepsilon)$}
\label{figHsquare}
\end{figure}

\begin{remark}\label{rem:collapses}{\em It is important to note that each square constituting $Q(\G,\varepsilon)$ corresponds to a 4-colored edge of $\Gamma$, while its sides are dual to $\{i,4\}$-colored cycles of $\G$ for each $i\in\Delta_3$ (see Figure \ref{figHsquare}).
Therefore any possible sequence of elementary collapses of the squares of $Q(\G,\varepsilon)$ can be encoded by an ordering of the corresponding 4-colored edges of $\G$. 
Note also that the collapse of the square corresponding to a 4-colored edge $e$ occurs if there exists  $i\in\Delta_3$ such that all squares having a side dual to the $\{4,i\}$-colored cycle containing $e$ have already collapsed, i.e if all 4-colored edges of this cycle precede $e$ in the ordering (see \cite{Casali-Cristofori gem-induced}, \cite{Casali-Cristofori trisection bis}).}
\end{remark}

\medskip

The existence of gem-induced trisections has been proved for a large class of manifolds, possibly comprehending all closed simply-connected ones.

\begin{proposition}\label{existence_gem-induced}  \  {\rm \cite{Casali-Cristofori trisection bis} } 
Let $M$ be a compact orientable $4$-manifold with empty or connected boundary.
If $M$ admits a handle-decomposition with no 3-handles, it also admits a gem-induced trisection. 
 
\noindent In particular, If $M$ is closed, it admits a gem-induced trisection if and only if it admits a handle decomposition lacking in 3- or 1-handles.
\end{proposition}

\begin{remark}{\em The above result has been proved by making use of an algorithm, presented in \cite{Casali-Cristofori Kirby-diagrams}, to construct a gem of a compact $4$-manifold represented by a Kirby diagram. The resulting graph turns out to give rise to a gem-induced trisection of the represented manifold. Moreover, its genus can also be expressed in terms of the combinatorial properties of the Kirby diagram.}
\end{remark}

In order to obtain trisections of a closed orientable $4$-manifold from any gem, when one of the 3-dimensional piece of the induced decomposition is not a handlebody, Martini and Toriumi (\cite{Martini-Toriumi}) recently suggested to perform  a kind of {\it stabilizations} along all edges of a fixed color of the gem (which is thought of as embedded in the triangulated manifold as the 1-skeleton of the dual complex). 
We can easily generalize this operation to the case of the manifold being non-orientable and/or with non-empty connected boundary.

\begin{definition}\label{def:stabilization}{\em Let $\G\in G_s^{(4)}$ be a gem of a compact 4-manifold with empty or connected boundary and $\e$ a cyclic permutation of $\Delta_4$.  If $\bar e$ is a $4$-colored\footnote{If $M$ is closed, edges of any color $c$ can be used, having care to exchange the roles of color $c$ and 4.} edge of $\G$ and $N_{\bar e}$ is a regular neighbourhood of $\bar e$ in $M$, the new decomposition $\mathcal  T^\prime(\Gamma, \varepsilon)=(H'_0,H'_1,H'_2)$ of $M$, where $H'_0=H_0\cup N_{\bar e}$ and $H'_i=\overline{H_i\setminus N_{\bar e}}$ for $i\in\{1,2\},$ is said to be obtained from  $\mathcal  T(\Gamma, \varepsilon)=(H_0,H_1,H_2)$ by {\it stabilization along} $\bar e.$   
\noindent Moreover, with little abuse of notation, $\mathcal  T^\prime(\Gamma, \varepsilon)$ itself will often be called the {\it stabilization along} $\bar e$  of $\mathcal  T(\Gamma, \varepsilon)$.}
\end{definition}

Obviously the above operation can be performed repeatedly along any number of different $4$-colored edges of $\G$; for simplicity, we will continue to denote the resulting decomposition by $(H'_0,H'_1,H'_2)$.
\medskip

By the above considerations and in view of the estimations of the trisection genus via gems, which will be discussed in the next section, it is useful to widen the definition of gem-induced trisection as follows.

\begin{definition}\label{def:geminduced}{\em A {\it  (generalized) gem-induced trisection} of a compact 4-manifold $M$ with empty or connected boundary is a decomposition $(H'_0,H'_1,H'_2)$ of $M$ obtained from $\mathcal  T(\Gamma, \varepsilon)$ ($\G\in G_s^{(4)}$ being a gem of $M$ and $\varepsilon$ a cyclic permutation of $\Delta_4$) by stabilizations along 4-colored edges and such that $H'_{12}=H'_1\cap H'_2$ is a 3-dimensional handlebody.  }
    \end{definition}

\begin{proposition}\label{prop:gem-induced} For each gem $\G\in G_s^{(4)}$ of a compact 4-manifold $M$ with empty or connected boundary and for each cyclic permutation $\varepsilon$ of $\Delta_4$, $\mathcal  T(\Gamma, \varepsilon)$ gives rise to a (generalized) gem-induced trisection of $M.$
\end{proposition}

\dimo Let $\mathcal  T(\Gamma, \varepsilon)=(H_0,H_1,H_2)$ be the decomposition of $M$ associated to a gem $\G\in G_s^{(4)}$ and a cyclic permutation $\e$ and let $(\widehat H_0,H_1,H_2)$ be the corresponding decomposition of $|K(\G)|=\widehat M$ obtained by capping off by a cone the ``outer'' boundary of $H_0$ (i.e. the boundary component of $\partial H_0$ not intersecting $H_1\cup H_2$). 

If $\mathcal  T(\Gamma, \varepsilon)$ is a gem-induced trisection (i.e. if $H_{12}=H_1 \cap H_2$ collapses to a graph), the statement trivially holds. Otherwise, let us see how  $\mathcal  T(\Gamma, \varepsilon)$ is modified by a stabilization along a 4-colored edge $\bar e$ of $\G.$   

Note that - since $\G$ is thought of as the 1-skeleton of the dual complex of $K(\G)$ - $\bar e$ is properly embedded in the complement of $\widehat H_0$ and also in the complement of $H_{12}.$
Therefore,  if $N_{\bar e}\cong\mathbb D^1\times\mathbb D^3$ is a (closed) regular neighbourhood of $\bar e$ in $K(\G)$, then $\partial N_{\bar e}\cong (\mathbb D^3\times\mathbb S^0)\cup (\mathbb S^2\times\mathbb D^1)$ intersects both $H_{01}=H_0\cap H_1$ and $H_{02}=H_0\cap H_2$ in a 3-ball. 
Moreover, $H_{12}$ intersects transversally the component of  $\partial N_{\bar e}$ homeomorphic to $\mathbb S^2\times\mathbb D^1$ and, therefore, splits it into two parts  both homeomorphic to $\mathbb D^1\times\mathbb D^2:$ the one intersecting $H_{01}$ (resp. $H_{02}$)  will be denoted by $h_1$ (resp. $h_2$).

Now, let us set $\widehat H'_0=\widehat H_0\cup N_{\bar e}$ and $H'_i=\overline{H_i\setminus N_{\bar e}}$ for $i\in\{1,2\}.$
Since $\widehat H_0\cap N_{\bar e} =\partial\widehat H_0\cap\partial N_{\bar e}\cong\mathbb D^3\times\mathbb S^0$, $\widehat H'_0$ is obtained from $\widehat H_0$ by attaching a 1-handle.

On the other hand, note that, for each $i\in\{1,2\}$,  $H_i\cap N_{\bar e}=h_i\times\mathbb D^1\cong\mathbb D^4$, while $\partial H_i\cap\partial N_{\bar e}\cong\mathbb D^3$, therefore $H_i$ can be seen as the boundary connected sum between $H'_i$ and a 4-ball. As a consequence $H'_i\cong H_i$. Furthermore, $H'_1$ (resp. $H'_2$) also collapses to the 1-dimensional subcomplex of $K(\G)$ generated by the $\e_0$- and $\e_2$-labeled (resp. $\e_1$- and $\e_3$-labeled) vertices.   

\begin{figure}   [h!]
     \centering
     \includegraphics[width=0.5\linewidth]{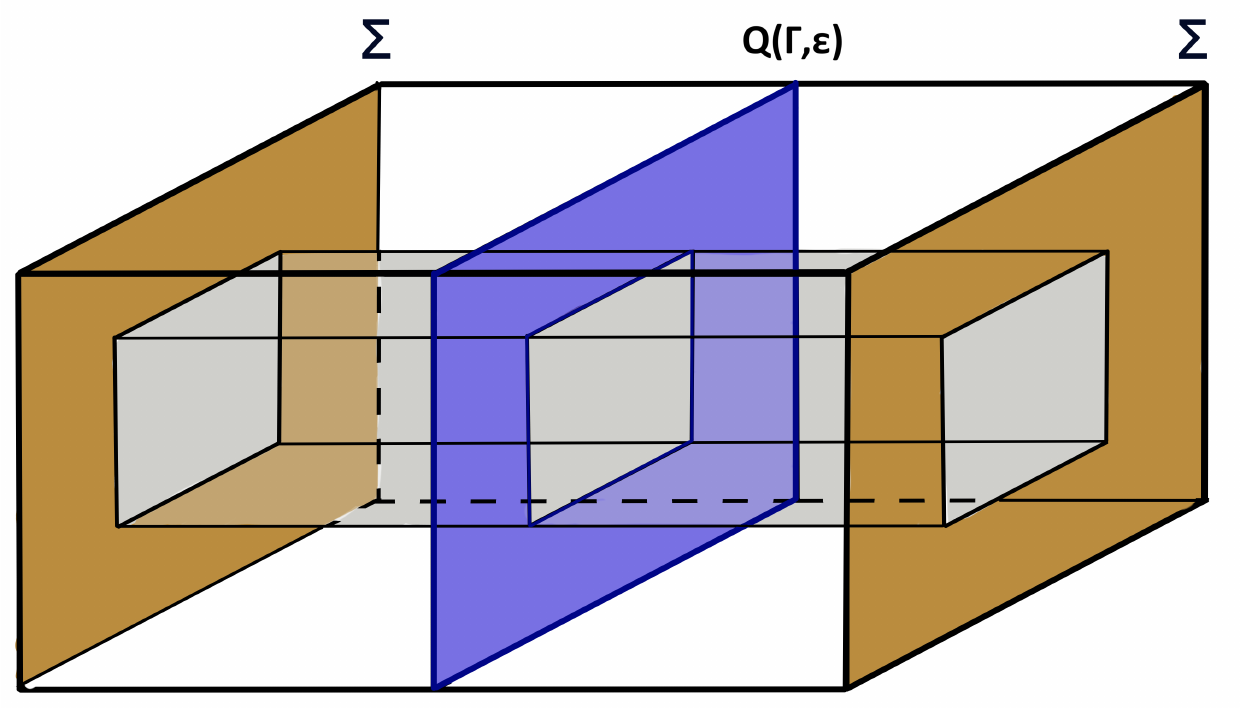}
     \caption{stabilization along a $4$-colored edge}
     \label{fig:carving_edge}
 \end{figure}

With regard to the 3-dimensional parts of the new decomposition $(\widehat H'_0,H'_1,H'_2)$ of $|K(\G)|$, the above observations imply that, for each $i\in\{1,2\}$, $H'_{0i}=\widehat H'_0\cap H'_i=(H_{0i}\setminus\partial N_{\bar e})\cup h_i$ is obtained by the attachment of one 1-handle to $H_{0i}.$ On the other hand, the carving of a 1-handle from $H_{12}$ yields $H'_{12}=H'_1\cap H'_2$  (see Figure \ref{fig:carving_edge}). The effect of this operation can be better understood by observing that the stabilization yields a ``hole'' in the square $q_{\bar e}$ of the spine $Q(\G,\e)$ of $H_{12}$ corresponding to the edge $\bar e$, thus ensuring the collapse of $q_{\bar e}$ to its 1-dimensional boundary.
 
Finally, it is not difficult to see that  $\Sigma'=\widehat H'_0\cap H'_1\cap H'_2$ is a closed connected surface resulting from the attachment of one 1-handle to $\Sigma=\widehat H_0\cap H_1\cap H_2$.

\smallskip

From this analysis it is clear that it certainly exists a sequence $(\bar e_1,\ldots , \bar e_k)$ 
such that, if we denote again by $(\widehat H'_0,H'_1,H'_2)$ the decomposition of $\widehat M$ obtained by repeatedly stabilizing $\mathcal T(\G,\e)$ along $(\bar e_1,\ldots , \bar e_k)$, then the spine $Q'(\G,\e)$ of $H'_{12}=H'_1\cap H'_2$ admits a sequence of elementary collapses to a 1-dimensional complex, i.e. $H'_{12}$ turns out to be a 3-dimensional handlebody (as well as $H'_{01}$ and $H'_{02}).$ 

By deleting from $\widehat H'_0$ a suitable regular neighbourhood of the (unique) 4-labelled vertex, a (generalized) gem-induced trisection $\mathcal T'(\G,\e)=(H'_0,H'_1,H'_2)$ of $M$ is finally  obtained.
\qed

\vskip-0.3cm

\vskip-0.5cm

\begin{remark} \label{rem:stabilization-dipole}{\em 
Given a 4-manifold $M$ with empty or connected boundary, a cyclic permutation $\e$ of $\Delta_4$ (with $e_4=4$) and a gem $\G\in G_s^{(4)}$ of $M$, let $\tilde \G$ be the 5-colored graph obtained from $\G$ by adding a {\it 4-dipole} on a 4-colored edge $\bar e$: $\tilde \G$ has two new vertices, joined by four edges of colors $0, 1, 2, 3$, and each of these new vertices is $4$-adjacent to one of the end-points of $\bar e$. $\tilde \G$ still represents $M$, since the associated pseudocomplex $K(\tilde \G)$ is obtained from $K(\G)$ by inserting a triangulated 4-ball.\footnote{See \cite[Section 4]{Ferri-Gagliardi-Grasselli} for the general notion of $h$-dipole ($0\le h \le n$) and related results.} 
We point out that a decomposition $\mathcal T\tilde\G,\e)=(\tilde H_0,\tilde H_1, \tilde H_2)$ of $M$ can still be constructed similarly to $\mathcal T(\G,\e)$, with the difference that $\tilde H_0$ as well as $\tilde H_{01}, \tilde H_{02}$ and $\tilde \Sigma$ all turn out to have two connected components, one of which is trivial (i.e. it is a disk or a sphere of the suitable dimension). 

Furthermore, the decomposition  $\mathcal T'(\G,\e)$ obtained by stabilization of $\mathcal T(\G,\e)$ along $\bar e$ can also be thought of as being obtained from $\mathcal T(\tilde \G,\e)$ by subsequently performing the operation of Definition \ref{def:stabilization} first along $\bar e'$ and then along $\bar e'',$ where $\bar e'$ and $\bar e''$ are the two 4-colored edges adjacent to the dipole. 

The operation along $\bar e'$ produces a boundary connected sum (resp. connected sum) between the different components of $\tilde H_0$ as well as of $\tilde H_{01}, \tilde H_{02}$ (resp. of $\tilde \Sigma$), while it carves a tunnel 
in $H_{12}$ between its boundary components.  On the other hand, the subsequent operation along $\bar e''$ is a stabilization according to Definition \ref{def:stabilization}, thus increasing the genus.

Furthermore, the $\{\e_0,\e_2\}$-cycle (resp. $\{\e_1,\e_3\}$-cycle) of the dipole is precisely the boundary of a compression disk for the 1-handle that is added to $H_{02}$ (resp. $H_{01}$) by the stabilization along $\bar e.$
}\end{remark} 
\medskip

\begin{remark} \label{rem:sequence}{\em Note that the construction in Proposition \ref{prop:gem-induced} of the resulting 
(generalized) gem-induced trisection is completely described by an ordering  $(\bar e_1,\ldots,\bar e_k, e_{k+1}, \ldots, e_p)$ of the 4-colored edges of $\G$, where $(\bar e_1,\ldots,\bar e_k)$ are the edges involved in the stabilizations, while, as already pointed out in Remark \ref{rem:collapses}, the collapse of $Q'(\G,\e)$ to a 1-dimensional complex exactly corresponds to an ordering of the remaining $4$-colored edges $(e_{k+1}, \ldots, e_p)$ with the property that, for each $j\in\{k+1,\ldots,p\}$, there exists $i\in\Delta_3$ such that all 4-colored edges of  the $\{4,i\}$-colored cycle containing $e_j$ belong to the set  $\{e_1,\ldots,e_j\}.$
}\end{remark}

\bigskip

\section{Estimations of the trisection genus via gems} \label{ss:proving-estimations} 

The generalization of the concept of gem-induced trisection
presented in Section 
\ref{s:gem-induced-trisections} 
naturally leads to generalize also the definition of the PL invariant {\it G-trisection genus} already introduced in \cite{Casali-Cristofori gem-induced}.

\begin{definition}\label{def_GT-genus} {\em 
Given a compact $4$-manifold $M$ with empty or connected  boundary, its  \emph{ (generalized) G-trisection genus} $g_{GT}(M)$ is defined as the minimum genus among all (generalized) gem-induced trisections of  $M.$
}
\end{definition}

In the closed case the G-trisection genus  is obviously an upper bound for the trisection genus and allows its estimation in two possible ways: either directly from gems of the manifold (possibly through stabilizations), or indirectly, as shown by the following result, which extends a similar one in \cite{Casali-Cristofori trisection bis}.

\begin{theorem}\label{trisection_from_gem-induced} \ 
\ Let \ $M$ \ be \ a \ compact  $4$-manifold whose boundary is a connected sum of sphere bundles over $\mathbb S^1$ and let $\bar M$  be its associated closed $4$-manifold. 
Then   any (generalized) gem-induced trisection of $M$ gives rise to a trisection of $\bar M$ with the same central surface.  

\noindent As a consequence, $$g_T(\bar M) \le g_{GT} (M).$$
\end{theorem}

\dimo Let $\G\in G_s^{(4)}$ be a gem of $M$ and $\e$  a cyclic permutation of $\Delta_4$. If  $\mathcal T(\G,\e)=(H_0,H_1,H_2)$ is already a gem-induced trisection of $M$, the result is Theorem 4.1 of \cite{Casali-Cristofori trisection bis}. 
Otherwise, let $\mathcal T'(\G,\e)=(H'_0,H'_1,H'_2)$ be a gem-induced trisection of $M$,  which is obtained by $k>0$ stabilizations from $\mathcal T(\G,\e)$; then $H'_0$ is obtained from $H_0$, i.e. a collar of $\partial M$,  by attaching $k$ 1-handles on its ``inner'' boundary component (the one intersecting $H_1\cup H_2$). Laudenbach-Poenaru's well-known result in \cite{Laudenbach-Poenaru}, or its non-orientable counterpart in \cite{Miller-Naylor}, ensures that the closed 4-manifold $\bar M$ is uniquely obtained by gluing a suitable 
handlebody to $\partial M$, i.e. to the ``outer'' boundary of $H'_0$; hence, a 4-dimensional handlebody $\bar H'_0$ is obtained from $H'_0.$ 
Therefore the statement is proved by noting that $(\bar H'_0,H'_1,H'_2)$ is a trisection of $\bar M$ whose central surface is $\bar H'_0\cap H'_1\cap H'_2=H'_0\cap H'_1\cap H'_2.$\qed

\medskip

A general result about the estimations of the trisection genus of any closed $4$-manifold via gems, as stated in Theorem \ref{intro:theorem} of the Introduction, can be obtained by analyzing the decompositions induced by gems with the aim of minimizing the number of stabilizations that are necessary to get a gem-induced trisection. From this perspective, the following result holds.

\begin{proposition}\label{GT-bounds}   \ 
Let $\G\in G_s^{(4)}$ be a gem of a compact $4$-manifold $M$ with empty or connected boundary and $\varepsilon$ a cyclic permutation of $\Delta_4$. Then $\mathcal  T(\Gamma, \varepsilon)$ gives rise to a (generalized) gem-induced trisection $\mathcal T^\prime (\Gamma, \varepsilon)$ of $M$ after stabilizations along a suitable number $k$ of 4-colored edges of $\G$, with $rk(\widehat M) \le k \le \rho_\e(\Gamma) - \rho_{\e}(\Gamma_{\hat 4})$, $rk(\widehat M)$ being the rank of $\pi_1(\widehat M)$.
\end{proposition}

\dimo  Given a cyclic permutation $\e$ and a gem $\G\in G_s^{(4)}$ of $M$, let $Q(\G,\e)$ be the 2-dimensional complex formed by the squares that are the bottom faces of the cubes constituting $H_{12} = H_1\cap H_2$ (see Figure \ref{figHsquare}). We are showing how to determine algorithmically a sequence of stabilizations of $\mathcal T(\G,\e)$ so as to yield a gem-induced trisection of $M$. 
As already pointed out in Remark \ref{rem:sequence}, any sequence of stabilizations and/or elementary collapses of $Q(\G,\e)$ can be represented by an ordering of the 4-colored edges of $\G$ corresponding to the involved squares of $Q(\G,\e).$

Let us first consider the (possibly disconnected) 3-colored graph $\G_{\{\e_0,\e_3,4\}}$. Since it is a 3-residue in a gem of an $n$-manifold, with $n\ge 3$, it represents the 2-sphere and therefore regularly embeds into 
$\mathbb S^2$ (\cite{Ferri-Gagliardi-Grasselli}); so the regions of the embedding are bounded by the $\{\e_0,\e_3\}$-, $\{\e_0,4\}$- and $\{\e_3,4\}$-colored cycles of $\G$.

Let now $\{\bar e_1,\ldots,\bar e_{\bar k}\}$ be a set of 4-colored edges of $\G$ given by the union of a spanning tree for each component of the graph obtained by shrinking to a point each $\{\e_0,\e_3\}$-colored cycle of $\G_{\{\e_0,\e_3,4\}}$. Note that $\bar k=g_{\e_0,4} - g_{\e_0,\e_3,4}.$

By the planarity of $\G_{\{\e_0,\e_3,4\}}$, all the remaining 4-colored edges of $\G$ can be ordered in a sequence $(e_{\bar k+1},\ldots,e_p)$ satisfying the condition of Remark \ref{rem:sequence} with $i\in\{\e_0,\e_3\}.$ 

Therefore, as already observed in the same remark, the ordering $(\bar e_1,\ldots,\bar e_{\bar k}, e_{\bar k+1},\ldots,e_p)$ encodes a collapsing sequence of $Q(\G,\e)$, which involves all squares of $Q(\G,\e)$. Hence the resulting decomposition of $M$ is a gem-induced trisection since all its 3-dimensional parts turn out to be handlebodies.

Furthermore, the combinatorial properties of 5-colored graphs representing compact 4-manifolds imply $\bar k = g_{\e_0,4} - g_{\e_0,\e_3,4} = \rho_\e(\G) -  \rho_{\e}(\Gamma_{\hat 4})$ (see Equation (5) of \cite[Proposition 3.1]{generalized-genus}), proving the upper bound for the number of stabilizations.

As regards the lower bound, let $\mathcal T'(\G,\e)=(H'_0,H'_1,H'_2)$ be a gem-induced trisection of $M$,  which is obtained by $k>0$ stabilizations from $\mathcal T(\G,\e)=(H_0,H_1,H_2)$, and let $(\widehat H'_0,H'_1,H'_2)$ be the associated decomposition of $\widehat M;$ we recall that $\widehat H'_0$ is obtained  by attaching $k$ 1-handles on a 4-ball (if $\partial M=\emptyset$) or on the cone over $\partial M$ (if $\partial M\neq\emptyset$), i.e. in both cases on a contractible space. Therefore $rk(\pi_1(\widehat H'_0))=k.$

On the other hand, since $\Sigma'=\widehat H'_0\cap H'_1\cap H'_2$ is a Heegaard surface both for $\partial H'_1$ and for $\partial H'_2$, the application of Seifert-Van Kampen's theorem first to the pair $(H'_1,H'_2)$ and then to the pair $(H'_1\cup H'_2,\widehat H'_0)$, proves that there is a surjection\footnote{The same argument is used in \cite[Remark 15]{Casali-Cristofori gem-induced} limited to the case $k=0$, but it can extended to the general one with few adjustments.} from $\pi_1(\widehat H'_0)$ to $\pi_1(\widehat M).$ Hence $rk(\widehat M)\leq k.$
\qed

\begin{corollary}\label{trisection_vs_regular-genus} \ 
For each compact $4$-manifold $M$ with non-empty connected boundary, then
$$\mathcal H(\partial M) + rk(\widehat M)\le g_{GT}(M) \le  \mathcal G(M).$$

\noindent where $\mathcal H(\partial M)$ is the Heegaard genus of $\partial M.$
\smallskip

\noindent 
For each closed $4$-manifold $M$, then
$$\mathcal G^{(\widehat \ )}(M) + rk(M)\le g_{GT}(M) \le  \mathcal G(M),$$
where \ $\mathcal G^{(\widehat \ )}(M) = min \{\rho_{\e}(\Gamma_{\hat i}) \ / \Gamma \ \text{gem of} \ M, \ \e \ \text{cyclic permutation of} \ \Delta_4, \ i \in \Delta_4\}.$
\end{corollary}

\dimo First, note that, if $M$ has non-empty boundary and $\mathcal T(\G,\e)=(H_0,H_1,H_2)$ is the decomposition associated to a gem $\G\in G_s^{(4)}$ and a cyclic permutation $\e$, then, as already pointed out in Section \ref{s:gem-induced-trisections}, $H_{01}$ and $H_{02}$ form a Heegaard splitting of $|K(\G_{\hat 4})|\cong\partial M$ of genus $\rho_{\varepsilon}(\G_{\hat 4})$; therefore  $\rho_{\e}(\Gamma_{\hat 4})\geq\mathcal H(\partial M).$ 
As a consequence, the lower and upper bounds for the G-trisection genus of $M$ directly come from Proposition \ref{GT-bounds}, since the genus of a gem-induced trisection obtained from $\mathcal T(\G,\e)$ by $k$ stabilizations is $\rho_{\varepsilon}(\G_{\hat 4})+k$. 

\smallskip
With regard to the closed case, it is sufficient to note that in Proposition \ref{GT-bounds} the role of color 4 can be taken by any other color.\qed 

From the above results it follows that, as already stated in the Introduction, the regular genus turns out to yield an estimation of the trisection genus. 

\bigskip

\noindent{\it Proof of Theorem \ref{intro:theorem}.\ } The claims are proved through the upper bounds of Corollary \ref{trisection_vs_regular-genus} and, in the case of boundary homeomorphic to a connected sum of (orientable or non-orientable) sphere bundles over $\mathbb S^1$, the subsequent application of Theorem \ref{trisection_from_gem-induced}.\qed

\bigskip

\section{From gems to trisection diagrams} \label{ss.from_gems_to_trisection_diagrams}
\par \noindent
 
Throughout this section, we will denote by $\#_m(\mathbb S^1 \otimes \mathbb S^2)$ the connected sum of $m$ (orientable or non-orientable)  $\mathbb S^2$-bundles over $\mathbb S^1$; moreover, in case $m=0,$ the same notation will indicate $\mathbb S^3$. 

\medskip 
According to  \cite{Miller-Naylor}, where both the orientable and non-orientable setting are 
taken into account, closed 4-manifolds can be  identified by means of suitable curves on surfaces, yielding trisections: 

\begin{definition} \label{def. trisection-diagram} {\em A $(g; k_0; k_1; k_2)$-{\it trisection diagram} is a 4-tuple $(\Sigma; \alpha, \beta,\gamma)$ such that the triples $(\Sigma; \alpha, \beta)$, $(\Sigma; \alpha, \gamma)$ and $(\Sigma; \beta,\gamma)$ are genus $g$ Heegaard diagrams for $\#_{k_0}(\mathbb S^1 \otimes \mathbb S^2)$, $\#_{k_1}(\mathbb S^1 \otimes \mathbb S^2)$ and $\#_{k_2}(\mathbb S^1 \otimes \mathbb S^2)$ respectively.}\footnote{Obviously, all $3$- and $4$-dimensional manifolds involved in Definition \ref{def. trisection-diagram}, as well as in the procedure to reconstruct the represented closed $4$-manifold, are orientable or not according to $\Sigma$.}   
\end{definition}  

In fact, such a diagram allows to uniquely determine  - up to diffeomorphism - a closed 4-manifold by the following steps:   
 \begin{itemize}
\item[-]  take the product $\Sigma \times \mathbb D^2$; 
  \item[-] take three genus $g$ 3-dimensional handlebodies $V_\alpha$, $V_\beta$ and $V_\gamma$, identified by the systems of curves $\alpha,$ $\beta$ and $\gamma$ respectively; 
  \item[-]  attach $V_\alpha \times I$ (resp. $V_\beta \times I)$ (resp. $V_\gamma \times I)$ to $\partial (\Sigma \times \mathbb D^2) = \Sigma \times \mathbb S^1$ along $\Sigma \times [-\epsilon, \epsilon]$ (resp. $\Sigma \times [\frac 2 3 \pi -\epsilon, \frac 2 3 \pi + \epsilon]$)  (resp. $\Sigma \times [\frac 4 3 \pi -\epsilon, \frac 4 3 \pi + \epsilon]$);  
  \item[-]  fill the three boundary components (isomorphic to $\#_{k_i}(\mathbb S^1 \otimes \mathbb S^2)$, for $i=0,1,2$) of the resulting compact $4$-manifold with three $4$-dimensional handlebodies (of genus $k_i$, for $i=0,1,2$).       
\end{itemize}

Note that the operation performed in the last step is well-defined due to 
the already cited theorem by Laudenbach and Poenaru (\cite{Laudenbach-Poenaru}), together with its non-orientable version proved in \cite{Miller-Naylor}: see also Remark \ref{rem:spine}.

\bigskip

In this section, we will prove how to get trisection diagrams directly from gems, via the procedures described in Sections \ref{s:gem-induced-trisections} and \ref{ss:proving-estimations}. 


Let us first consider a compact $4$-manifold $M$ that is either closed or  whose boundary is a connected sum of sphere bundles over $\mathbb S^1$. In this case, Proposition \ref{prop:gem-induced} ensures that from any gem $\Gamma\in G_s^{(4)}$ of $M$ and any cyclic permutation $\e$ of $\Delta_4$, a gem-induced trisection $\mathcal T^\prime(\G, \e)$ of $M$ can be always obtained by applying $k \ge 0$ stabilizations\footnote{
Note that,
if $k=0$, the decomposition $\mathcal T(\G, \e)$ itself is a gem-induced trisection, and hence $\mathcal T^\prime(\G, \e)=\mathcal T(\G, \e)$.} to the decomposition $\mathcal T(\G, \e)$; moreover, by  Remark \ref{rem:sequence},  $\mathcal T^\prime(\G, \e)$ can be encoded by an ordering of all 4-colored edges of $\G$ 
$$(\bar e_1,\ldots,\bar e_k, e_{k+1}, \ldots, e_p),$$  with the property that, for each $j\in\{k+1,\ldots,p\}$, 
$$\begin{aligned} \text{there exists } i\in\Delta_3 \text{ such that} & \ \text{all 4-colored edges of  the }  \\ \{4,i\}\text{-colored cycle} \ C_j \ \text{containing } e_j  & \ \text{belong to the set }\{e_1,\ldots,e_j\}.\end{aligned}$$

Now, let $F_\e(\Gamma_{\hat 4})^{(k)}$ be the surface obtained from the embedding surface $F_\e(\Gamma_{\hat 4})$ of the subgraph $\Gamma_{\hat 4}$ (see Proposition \ref{reg_emb}) by adding $k\ge 0$ 1-handles along the edges $\bar e_1,\ldots,\bar e_k$ of the above ordering and let $K_{\e_1 \e_3}$ (resp. $K_{\e_0 \e_2}$) be the 1-dimensional subcomplex of $K(\G_{\hat 4})$ generated by the vertices labeled $\e_1$ and $\e_3$ (resp. $\e_0$ and $\e_2$).

\medskip
The following statement yields a trisection diagram which identifies either the trisection $\mathcal T^\prime(\G, \e)$ of $M$ (if $M$ is closed) or the trisection $\bar{\mathcal T}^\prime(\G, \e)$ of the associated closed 4-manifold $\bar M$ (if $\partial M$ is a connected sum of sphere bundles over $\mathbb S^1$).   

\begin{proposition} \label{trisection_diagrams_gem-induced(chiuse)} \  
Given a compact $4$-manifold $M$ 
such that $\partial M \cong \#_m(\mathbb S^1 \otimes \mathbb S^2)$, $m \ge 0$, let $\ \Gamma\in G_s^{(4)}$ be a gem of $M$ and $\varepsilon=(\e_0,\e_1, \e_2, \e_3, 4)$ a cyclic permutation of $\Delta_4.$ Furthermore, if  $(\bar e_1,\ldots,\bar e_k, e_{k+1}, \ldots, e_p)$ is an ordering of the 4-colored edges of $\G$ encoding a gem-induced trisection of $M,$ let $\tilde \G$ be the 5-colored graph obtained from $\G$ by adding a 4-dipole on $\bar e_i$, for each $i\in\{1,\ldots,k\}.$

Then, with the above notations, a trisection diagram for the associated closed $4$-manifold $\bar M$ is given by the following three systems of curves $\{ \alpha, \beta, \gamma\}$ on $F_\e(\Gamma_{\hat 4})^{(k)}:$  
\begin{itemize}
\item[-] the $\alpha$ curves (resp. $\beta$ curves) consist of all $\{\e_0,\e_2\}$-colored cycles (resp. $\{\e_1,\e_3\}$-colored cycles) of $\tilde \G$, 
but those corresponding to a spanning tree of $K_{\e_1 \e_3}$ (resp. $K_{\e_0 \e_2}$);
%
\item[-] the $\gamma$ curves form a subset of the $\{4,i\}$-colored cycles of $\tilde \Gamma$ ($i\in \Delta_3$)  which can be 
obtained by first deleting the $\{4,i\}$-colored cycles of $\tilde \G$ arising from the $C_j$'s, 
then by also deleting the cycles of $\tilde \G$ 
corresponding to a spanning tree of $Q_1 \setminus \{c_{k+1}, \dots, c_p\},$ where $Q_1$ denotes the $1$-skeleton of $Q(\G,\e)$ and $c_j$ the edge of $Q(\G,\e)$ corresponding to the cycle $C_j$. 
\end{itemize} 
\end{proposition}

\dimo 
First of all, let us note that the $\{4,i\}$-cycles of $\tilde \G$ correspond bijectively to those of $\G$ (actually they simply differ by some $4$-colored edges being split by $i$-colored edges of the added dipoles): for this reason, for sake of conciseness, in the following we will talk indifferently of $\{4,i\}$-cycles without any specification of the gem they belong to. On the other hand, $\tilde \G$ has $k$ $\{\e_0,\e_2\}$-cycles (resp. $\{\e_1,\e_3\}$-cycles) more than $\G$: they have length two and are contained in the $k$ added 4-dipoles. 
 
\smallskip

Let us recall that $\mathcal T(\Gamma, \varepsilon)=(H_{0},H_{1},H_{2})$, where $H_{1}$ and  $H_{2}$ are $4$-dimensional handlebodies, while  $H_{0}$ is either the $4$-disk obtained by coning over the triangulation $K(\Gamma_{\hat 4})$ of the $3$-sphere (if $m=0$) or a collar of the triangulation $K(\Gamma_{\hat 4})$ of $\partial M \cong \#_m(\mathbb S^1 \otimes \mathbb S^2)$ (if $m>0$).
Moreover, as pointed out by Martini and Toriumi in the closed orientable case\footnote{See \cite[Paragraph 4.6]{Martini-Toriumi}, and in particular the caption of Figure 17.}, the decomposition $\mathcal T(\Gamma, \varepsilon)$ of $M$ induces a quadrangulation of the surface $\Sigma=H_0\cap H_1 \cap H_2$ (given by the union of the squares depicted in orange in Figure \ref{fig:trisec_pieces}), which is dual to the cellular subdivision of the embedding surface $F_\e(\Gamma_{\hat 4})$ given by $\Gamma_{\hat 4}$ itself; hence, $\Sigma$ and  $F_\e(\Gamma_{\hat 4})$ can be identified.

Now, if $\mathcal  T^\prime(\Gamma, \varepsilon) = (H^\prime_{0},H^\prime_{1},H^\prime_{2})$ is the gem-induced trisection of $M$ obtained from $\mathcal T(\Gamma, \varepsilon)$ by stabilization along the edges $(\bar e_1, \dots, \bar e_k)$ \ ($k \ge 0$), the identification between $\Sigma^\prime= H^\prime_{0} \cap H^\prime_{1} \cap H^\prime_{2}$ and $F_\e(\Gamma_{\hat 4})^{(k)}$ directly follows by construction (see the proof of Proposition \ref{prop:gem-induced}).

Hence, in order to obtain a trisection diagram which identifies the trisection $\bar{\mathcal  T}^\prime(\Gamma, \varepsilon) = (\bar{H}^\prime_{0},H^\prime_{1},H^\prime_{2})$ of the associated closed $4$-manifold $\bar M$ (which is exactly the trisection $\mathcal  T^\prime(\Gamma, \varepsilon)$ of $M$, if $M$ is closed), whose central surface is the same as $\mathcal  T^\prime(\Gamma, \varepsilon)$, we have to determine on $F_\e(\Gamma_{\hat 4})^{(k)}$ 
a system of independent curves realizing the attachments of the 3-dimensional handlebodies  $\bar{H}^\prime_{0} \cap H^\prime_{1} =  H^\prime_{0} \cap H^\prime_{1}=H^\prime_{01}$, $\bar{H}^\prime_{0} \cap H^\prime_{2} =  H^\prime_{0} \cap H^\prime_{2}=H^\prime_{02}$ and $H^\prime_{1} \cap H^\prime_{2}=H^\prime_{12}$ .    

\medskip

As already pointed out, the union of the red (resp. green) edges to which the prisms constituting $H_{02}=H_0 \cap H_2$ (resp. $H_{01}=H_0 \cap H_1$) collapse (see Figure \ref{fig:trisec_pieces}) is a 1-dimensional spine of the 3-dimensional handlebody $H_{02}$ (resp. $H_{01}$); moreover, it is not difficult to check that this spine is isomorphic to $K_{\e_1 \e_3}$ (resp. $K_{\e_0 \e_2}$). As a consequence, a ``minimal" spine of $H_{02}$ (resp. $H_{01}$) is obtained by shrinking to a point the red (resp. green) edges corresponding to a spanning tree of $K_{\e_1 \e_3}$ (resp. $K_{\e_0 \e_2}$), and a system of (independent) meridian curves of $H_{02}$ (resp. $H_{01}$) consists of curves which are dual to the edges of that minimal spine. 
In this regard, note that each $\{\e_0,\e_2\}$-colored cycle (resp. $\{\e_1,\e_3\}$-colored cycle) of $\G$, is dual in $K(\Gamma_{\hat 4})$ to one of the above red (resp. green) edges (see \cite[Figure 17]{Martini-Toriumi}) and it is embedded in $F_\e(\Gamma_{\hat 4})$ (and then also in $F_\e(\Gamma_{\hat 4})^{(k)}$). 

Moreover, as pointed out in Remark \ref{rem:stabilization-dipole}, for each $i\in\{1,\ldots,k\}$, the length two $\{\e_0,\e_2\}$-cycle (resp. $\{\e_1,\e_3\}$-cycle) of $\tilde \G$ belonging to the 4-dipole added on $\bar e_i$ yields a meridian curve for the 1-handle of $H^\prime_{02}$ (resp. $H^\prime_{01}$) arising from the stabilization along $\bar e _i$ and it is obviously embedded in $F_\e(\Gamma_{\hat 4})^{(k)}$, too. 
The statement regarding $\alpha$ (resp. $\beta$) curves now easily follows.


\medskip 
As regards $H^\prime_{12}$, it is already known - as pointed out in Remark \ref{rem:sequence} - that its 2-dimensional spine $Q^\prime(\G, \e)$ (obtained from the spine $Q(\G, \e)$ of $H_{12}=H_1\cap H_2$ by performing a ``hole" in each square $q_{\bar e_j}$, for each $j=1, \dots, k$) collapses to a 1-dimensional subcomplex: more precisely, all squares $q_{\bar e_j},$ for $j=1, \dots, k$, collapse to their boundaries, while each square $q_{e_j}$, for $j=k+1, \dots, p$, collapses from its free edge $c_j$, corresponding to the $\{4,i\}$-colored cycle $C_j$ containing $e_j$ so that all its $4$-colored edges belong to the set $\{e_1,\ldots,e_j\}.$ 
Hence, a 1-dimensional spine of $H^\prime_{12}$ is obtained from the $1$-skeleton of $Q(\G, \e)$, by deleting all edges $c_j$, for $j=k+1, \dots, p$. 
Moreover, a ``minimal" 1-dimensional spine of $H^\prime_{12}$ can be obtained by shrinking to a point the edges belonging to a spanning tree.  As a consequence, a system of (independent) meridian curves, realizing the attachment of $H^\prime_{12},$ consists of the set $X$ of curves dual to the edges of the above minimal spine.  

The statement regarding the $\gamma$ curves follows by recalling that each element of $X$ corresponds to a $\{4,i\}$-cycle and by proving that it can be projected on $F_\e(\Gamma_{\hat 4})^{(k)}$. In fact, all $i$-colored edges ($i \in \Delta_3$) obviously embed into $F_\e(\Gamma_{\hat 4})$, while each $\bar e_j$ ($j\in\{1,\ldots,k\}$) embeds into $F_\e(\Gamma_{\hat 4})^{(k)}$ by construction; in particular, this embedding can be better understood by means of the added dipoles, since the stabilization along $\bar e_j$ can be factorized as described in Remark \ref{rem:stabilization-dipole}. Moreover, the given ordering of the 4-colored edges enables to project on $F_\e(\Gamma_{\hat 4})^{(k)}$ also the 4-colored edge $e_j$, for each $j=k+1, \ldots, p,$ via the bicolored cycle $C_j$ containing it
\footnote{Note that the argument is consistent with both \cite[Section 4.4]{Martini-Toriumi}, where $k=p$ is always chosen, and \cite[Proposition 3.20]{CC-Contemporary}, where $k=0$ holds.}.
%
%
\qed 

\begin{figure} [H] 
    \centering
    \includegraphics[width=0.6\linewidth]{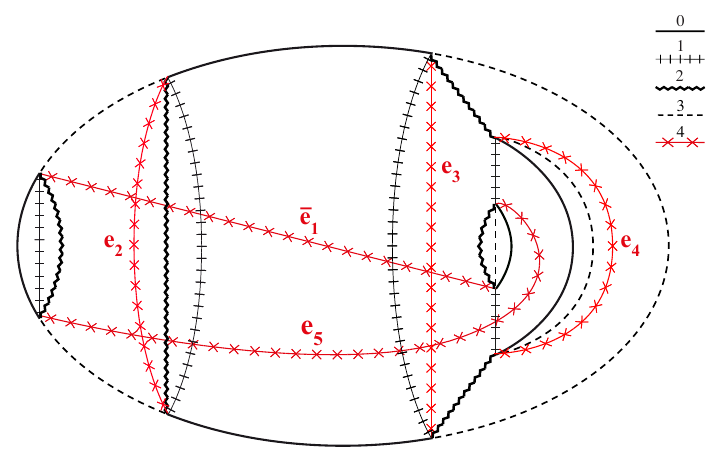}
   \caption{The standard order ten crystallization of $\mathbb S^3 \times \mathbb S^1$,  with an ordering of $4$-colored edges giving rise to a (generalized) gem-induced trisection}   
   \label{fig:S3xS1_sequence}
\end{figure}

\begin{figure} [H]
    \centering
    \includegraphics[width=0.3\linewidth]{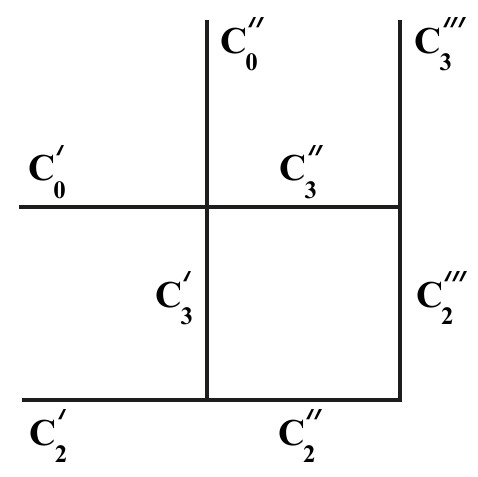}
   \caption{The 1-dimensional complex arising from $Q(\Gamma, \varepsilon)$ after stabilization and collapses}
   \label{fig:K1_S3xS1}
\end{figure}

\begin{figure} [H] 
    \centering
    \includegraphics[width=0.6\linewidth]{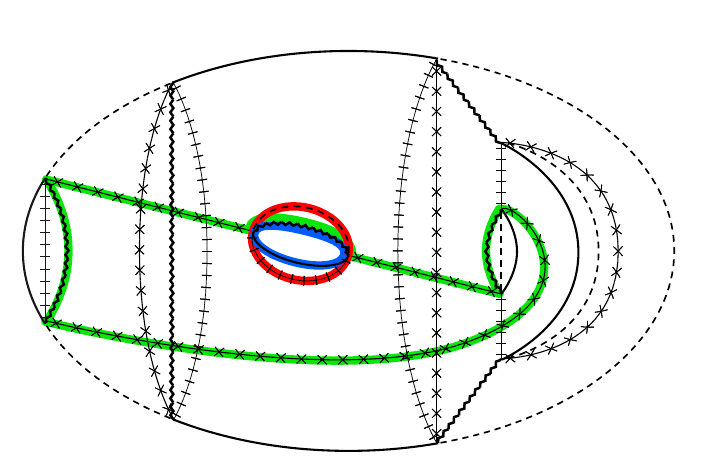}
    \caption{The curves of the trisection diagram  of $\mathbb S^3 \times \mathbb S^1$ as cycles of the gem $\tilde \G$} 
   \label{fig:S3xS1_trisection-diagram(highlighted)}
\end{figure}

\begin{figure} [H] 
\centering
\begin{minipage}[c]{0.4\textwidth}
\hspace*{0.5truecm}
\includegraphics[width=1.2\linewidth]{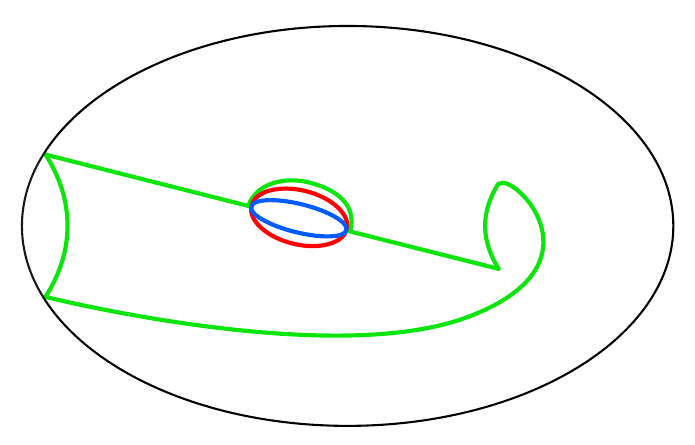}
\end{minipage}
\hfill
\begin{minipage}[c]{0.1\textwidth}
\centering
\hspace*{2truecm}
\Large $\approx$
\end{minipage}
\hfill
\begin{minipage}[c]{0.4\textwidth}
\centering
\includegraphics[width=0.7\linewidth]{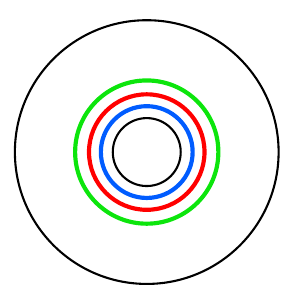}\end{minipage}
\hspace*{0.5truecm}
\caption{The associated trisection diagram  of $\mathbb S^3 \times \mathbb S^1$}
\label{fig:S3xS1_standard_diagram}
\end{figure}

\begin{example} \label{ex:diagram_S3xS1} {\em
Figures \ref{fig:S3xS1_sequence},  \ref{fig:S3xS1_trisection-diagram(highlighted)} and \ref{fig:S3xS1_standard_diagram}  
show the application of 
Proposition \ref{trisection_diagrams_gem-induced(chiuse)} 
to the gem $\Gamma$ of Figure \ref{fig:S3xS1_sequence}, which is the standard order ten crystallization of $\mathbb S^3 \times \mathbb S^1$. $\G$ is proved to admit a (generalized) gem-induced trisection - according to Definition \ref{def:geminduced} -, with respect to $\e=(0,1,2,3,4)$ and the depicted ordering of its $4$-colored edges, $(\bar e_1, e_2, e_3, e_4, e_5)$. In fact, if a stabilization along the edge $\bar e_1$ is performed, all squares of $Q(\Gamma, \varepsilon)$ 
subsequently collapse from their free edges corresponding 
to $\{c_i,4\}$-colored cycles 
with $c_2=c_3=1$, $c_4=0$, $c_5=1$.
The result is a 1-dimensional complex consisting of eight edges:
\begin{itemize}
    \item [-] three corresponding to all $\{2,4\}$-colored cycles of $\Gamma$: $C'_{2},\ C''_2$ and $C'''_2$, say, containing  $e_2$,$\ \{\bar e_1, e_5\}$ and  $\{e_3, e_4\}$ respectively;
    \item [-] three corresponding to all $\{3,4\}$-colored cycles of $\Gamma$: $C'_3,\ C''_3$ and $C'''_3$, say, containing $\{\bar e_1, e_2, e_5\},\ e_3$ and $e_4$ respectively;
     \item [-] two corresponding to the $\{0,4\}$-colored cycles $C'_{0}$ and $C''_{0}$, say, containing $\{e_2, e_3\}$ and $\{\bar e_1, e_5\}$ respectively.    
\end{itemize}

For a picture of the above complex see Figure \ref{fig:K1_S3xS1}, where each edge is identified by the name of the corresponding bicolored cycle.

By shrinking to a point a maximal tree, 
a unique edge remains, for example the one corresponding to the cycle $C''_2$, which is highlighted in green in Figure \ref{fig:S3xS1_trisection-diagram(highlighted)} as cycle of the gem $\tilde \G$ obtained from $\G$ by adding a 4-dipole  on $\bar e_1$. Hence, $C''_2$ identifies the (unique) curve of the system $\gamma$ in the associated trisection diagram of genus
one, whose surface is obtained, via the stabilization, from the 2-sphere where $\Gamma_{\hat 4}$ regularly embeds.
On the other hand, the systems $\alpha$ and $\beta$ consist of the blue and red curves determined by the $\{0,2\}$- and $\{1,3\}$-cycle of the dipole respectively.   


It is easy to check that the obtained trisection diagram exactly coincides with the standard trisection diagram of $\mathbb S^3 \times \mathbb S^1$  (see Figure \ref{fig:S3xS1_standard_diagram}). 
}
\end{example}

\medskip 

\medskip 

\begin{figure} [h] 
    \centering
    \includegraphics[width=0.75\linewidth]{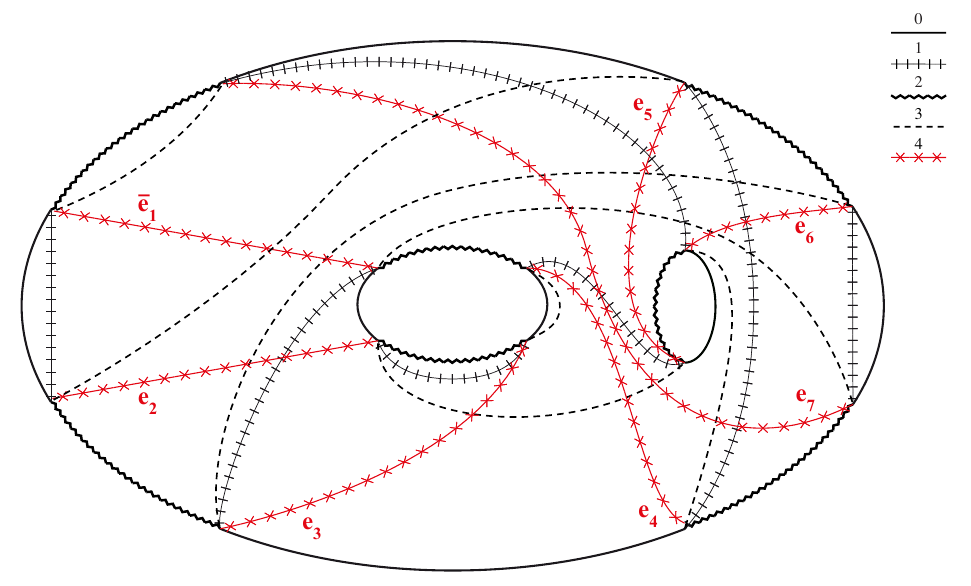}
    \caption{A gem of $\mathbb S^2 \times \mathbb S^2$, with an ordering of $4$-colored edges giving rise to a (generalized) gem-induced trisection}
   \label{fig:S2xS2_dabucare_grafo}
\end{figure}

\begin{figure} [h] 
    \centering
    \includegraphics[width=0.75\linewidth]{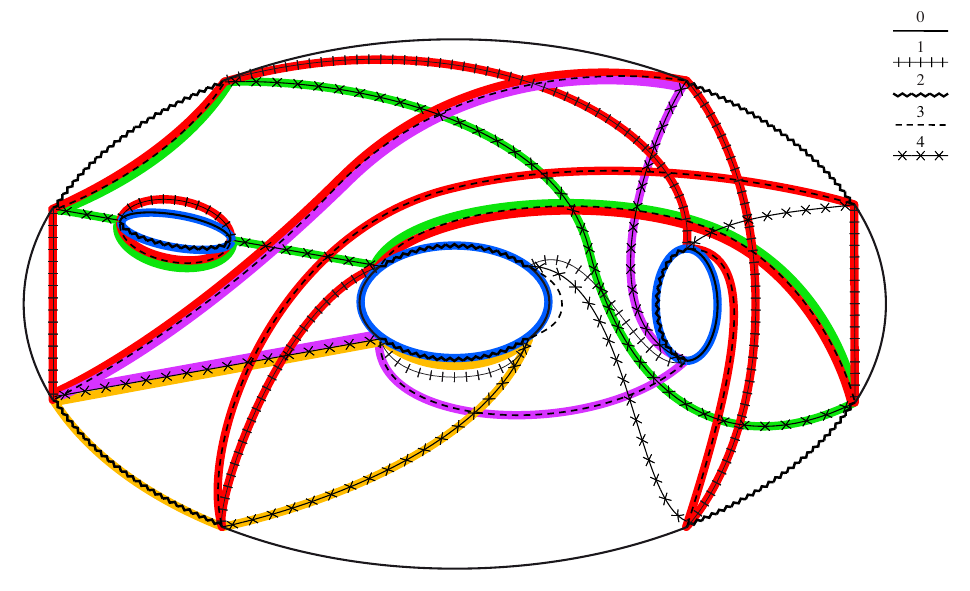}
    \caption{The curves of the trisection diagram  of $\mathbb S^2 \times \mathbb S^2$ as cycles of the gem $\tilde \G$} 
   \label{fig:S2xS2_dabucare_cicli-evidenziati}
\end{figure}

\begin{figure} [H] 
    \centering
    \includegraphics[width=0.75\linewidth]{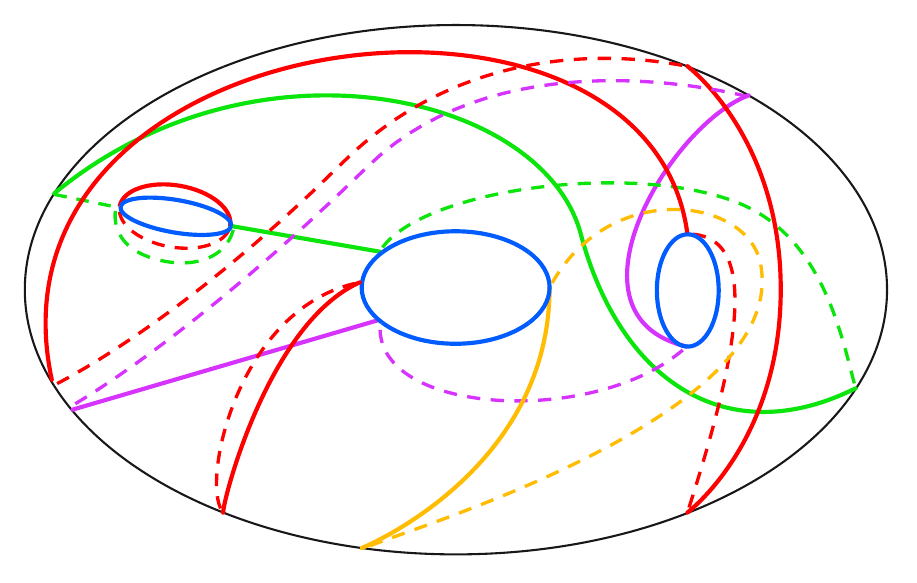}
    \caption{The associated trisection diagram  of  $\mathbb S^2 \times \mathbb S^2$}
   \label{fig:S2xS2_dabucare_trisection-diagram}
\end{figure}

\begin{example} {\rm \ \\
Figure \ref{fig:S2xS2_dabucare_grafo} shows a crystallization $\G$ of order 14 of $\mathbb{S}^2 \times \mathbb{S}^2$ that requires one stabilization to induce a trisection with respect to the permutation $\e=(0,1,2,3,4)$. For example, after stabilizing along the 4-colored edge $\bar e_1$, the depicted ordering of the remaining 4-colored edges yields a collapse of $Q(\Gamma, \varepsilon)$ to a 1-dimensional complex. After shrinking a maximal tree, three edges remain, whose corresponding bicolored cycles are, for instance, those highlighted in yellow, fuchsia 
 and green in Figure \ref{fig:S2xS2_dabucare_cicli-evidenziati} as cycles of the gem $\bar\G$ obtained by adding a 4-dipole on $\bar e_1$. Hence these three curves form the system $\gamma$, while the curves highlighted in blue and red constitute the systems $\alpha$ and $\beta$ respectively. Figure \ref{fig:S2xS2_dabucare_trisection-diagram} shows the embedding of the curves of the three systems in the (genus 3) orientable surface obtained after the stabilization from the surface where $\G_{\hat 4}$ regularly embeds with respect to $\e$, and hence it yields a trisection diagram of genus 3 of $\mathbb{S}^2 \times \mathbb{S}^2.$
}
\end{example}

\begin{figure} [h] 
    \centering
    \includegraphics[width=0.55\linewidth]{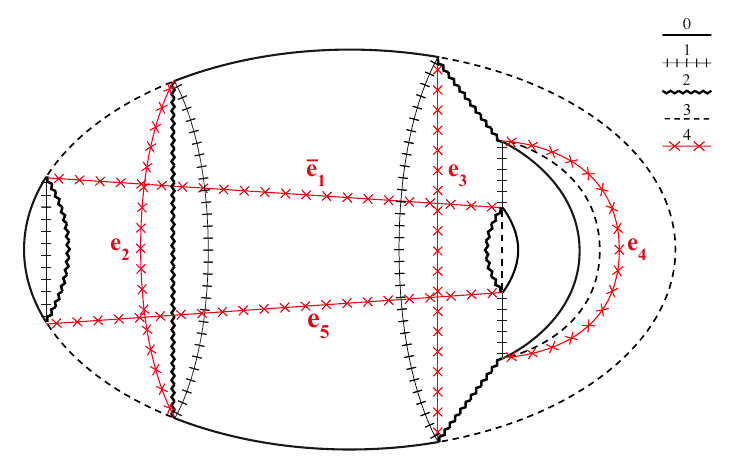}
    \caption{The standard order ten crystallization of $\mathbb S^3\tilde\times \mathbb S^1$, with an ordering of $4$-colored edges giving rise to a (generalized) gem-induced trisection}
   \label{fig:S3xS1_nonor}
\end{figure}

\begin{example} {\rm \ \\
In full analogy with Example \ref{ex:diagram_S3xS1}, 
Proposition \ref{trisection_diagrams_gem-induced(chiuse)} can be applied to the standard order ten crystallization of $\mathbb S^3 \tilde \times \mathbb S^1$ depicted in Figure \ref{fig:S3xS1_nonor}: the only difference is that the stabilization involves the attachment of a non-orientable handle (since the end-points of the edge $\bar e_1$ belong to the same bipartition class of  $\Gamma_{\hat 4}$), so the trisection surface is a non-orientable genus one surface. It is easy to check that the obtained trisection diagram exactly coincides with the standard trisection diagram of $\mathbb S^3 \tilde \times \mathbb S^1$ (see \cite{Miller-Naylor}). 
}
\end{example}

As already pointed out in \cite{CC-Contemporary}, the particular extension to the boundary case of the notion of trisection 
given by gem-induced trisections, suggests also a possible extension of the notion of trisection diagram to simply-connected compact $4$-manifolds with connected boundary. 

\begin{definition} \label{def. G-trisection-diagram}  \ \ {\rm (\cite{CC-Contemporary})} \ \ 
 A \emph{G-trisection diagram}  of genus $g$ is a 4-tuple $(\Sigma; \alpha, \beta, \gamma)$, where $\Sigma$ is a genus $g$ orientable surface and $\alpha$, $\beta$, $\gamma$ are complete systems of meridian curves on $\Sigma$, such that: 
\begin{itemize}
    \item[-] 
  $(\Sigma; \alpha, \gamma)$ and $(\Sigma; \beta, \gamma)$ are (genus $g$) Heegaard diagrams of $\#_{k_1}(\mathbb S^1 \times \mathbb S^2)$ and $\#_{k_2}(\mathbb S^1 \times \mathbb S^2)$ respectively; 
  \item[-] 
 $(\Sigma; \alpha, \beta)$ is a (genus $g$) Heegaard diagram of a closed connected 3-manifold.  
\end{itemize}

\end{definition}

In \cite{CC-Contemporary}, the above extension has been applied only in the case of manifolds directly admitting gem-induced trisections. Here, by taking advantage of the results of Sections \ref{s:gem-induced-trisections} and \ref{ss:proving-estimations}, the whole class of simply-connected compact $4$-manifolds with connected boundary can be considered. 

 \begin{proposition}   \label{trisection_diagrams_simply-connected_bounded}
Let $M$ be a simply-connected compact $4$-manifold with connected boundary. 
Then, a G-trisection diagram $(\Sigma; \alpha, \beta, \gamma)$ of genus $g$ exists, which uniquely identifies $M$, such that $\mathcal H(\partial M) + rk(\widehat M) \le g \le \mathcal G(M)$ and $(\Sigma; \alpha, \beta)$ is a Heegaard diagram of  $\#_k(\mathbb S^1\times \mathbb S^2)\#\partial M$ ($k \ge 0)$.     
\end{proposition} 


\dimo 
Let $\Gamma\in G_s^{(4)}$ be a gem of $M$ that realizes its regular genus, i.e. $\rho_\e(\G)=\mathcal G(M)$, \ $\e=(\e_0,\e_1, \e_2, \e_3, 4)$ being a suitable cyclic permutation of the color set $\Delta_4$. 
According to Proposition \ref{prop:gem-induced} and Proposition \ref{GT-bounds}, let ${\mathcal T}^\prime(\Gamma, \e)= (H^\prime_0, H^\prime_1, H^\prime_2)$ be the (generalized) gem-induced trisection of $M$ obtained from the subdivision $\mathcal T(\Gamma, \e)= (H_{0},H_{1},H_{2})$ via $k$ stabilizations, with $rk(\widehat M) \le k \le \rho_\e(\G) - \rho_\e(\G_{\hat 4})$. Its genus, equal to $\rho_\e(\G_{\hat 4}) + k$, obviously satisfies both bounds of the statement. 

By construction,  if $\alpha, \beta, \gamma$ are attaching curves on $\Sigma^\prime= H^\prime_0 \cap H^\prime_1 \cap H^\prime_2$ of $H^\prime_{01}=H^\prime_0 \cap H^\prime_1$, $H^\prime_{02} = H^\prime_0 \cap H^\prime_2$ and $H^\prime_{12}=H^\prime_1 \cap H^\prime_2$ respectively, then $(\Sigma^\prime; \alpha, \beta, \gamma)$ turns out to be a G-trisection diagram (see \cite[Remark 15]{Casali-Cristofori gem-induced} for details). 

\smallskip
On the other hand, 
let $W$ be the compact 4-manifold obtained by attaching a collar of two genus $g$ 3-dimensional handlebodies (corresponding to $H^\prime_{01}$ and $H^\prime_{02}$) to $\Sigma^\prime \times \mathbb D^2$  according to the curves $\alpha$ and $\beta$ respectively, in full analogy of the procedure described after Definition \ref{def. trisection-diagram}.  Since $M$ can be re-obtained from $W$ by adding a suitable number of 
3-handles (corresponding to $H^\prime_1$ and $H^\prime_2$, as well as to the $k$ 1-handles added to the ``inner" boundary of ${H}_{0}\cong\partial M\times I$ 
to obtain ${H}^\prime_{0}$), and $M$ is simply-connected with connected boundary by hypothesis, Theorem 1 in \cite{[Trace 1982]}  ensures that $W$ (or, equivalently, the G-trisection diagram) uniquely determines $M$.
\vskip-0.7cm\ \qed 

The following result can be proved by the same arguments as for 
Proposition \ref{trisection_diagrams_gem-induced(chiuse)}, with slight modifications.

\begin{proposition} \label{trisection_diagrams_gem-induced(bordo)} \  
Given a simply-connected compact $4$-manifold $M$ with connected boundary,  let $\Gamma\in G_s^{(4)}$ be a gem of $M$ and $\varepsilon=(\e_0,\e_1, \e_2, \e_3, 4)$ a cyclic permutation of $\Delta_4.$ Furthermore, if $(\bar e_1,\ldots,\bar e_k, e_{k+1}, \ldots, e_p)$ is an ordering of the 4-colored edges of $\G$ encoding a gem-induced trisection of $M,$ let $\tilde \G$ be the 5-colored graph obtained from $\G$ by adding a 4-dipole on $\bar e_i$, for each $i\in\{1,\ldots,k\}.$

Then, a G-trisection diagram for $M$ is given by the following three systems of curves  $\{ \alpha, \beta, \gamma\}$ on $F_\e(\Gamma_{\hat 4})^{(k)}:$  
\begin{itemize}
\item[-] the $\alpha$ curves (resp. $\beta$ curves) consist of all $\{\e_0,\e_2\}$-colored cycles (resp. $\{\e_1,\e_3\}$-colored cycles) of $\tilde \G$, 
but those corresponding to a spanning tree of $K_{\e_1 \e_3}$ (resp. $K_{\e_0 \e_2}$);
%
\item[-] the $\gamma$ curves form a subset of the $\{4,i\}$-colored cycles of $\tilde \Gamma$ ($i\in \Delta_3$)  which can be 
obtained by first deleting the $\{4,i\}$-colored cycles of $\tilde \G$ arising from the $C_j$'s, 
then by also deleting the cycles of $\tilde \G$ 
corresponding to a spanning tree of $Q_1 \setminus \{c_{k+1}, \dots, c_p\},$ where $Q_1$ denotes the $1$-skeleton of $Q(\G,\e)$ and $c_j$ the edge of $Q(\G,\e)$ corresponding to the cycle $C_j$. 
\end{itemize} 
\end{proposition} 
\vskip-0.8cm  \ \qed 

\begin{remark}
   {\rm  Proposition \ref{trisection_diagrams_gem-induced(chiuse)}  and Proposition \ref{trisection_diagrams_gem-induced(bordo)} generalize results presented in \cite{CC-Contemporary} for gems that directly yield 
   gem-induced trisections: see \cite[Propositions 3.20 and 3.25]{CC-Contemporary}.  
In particular, in \cite[Propositions 3.22 and 3.26]{CC-Contemporary}, further  simplifications are obtained for the special case of gems arising from framed links, by taking advantage of the algorithmic procedure described in \cite{Casali-Cristofori Kirby-diagrams}.
}
\end{remark}

\bigskip \bigskip 

\noindent {\bf Acknowledgements:\ }  This work was supported by GNSAGA of INDAM and by the University of Modena and Reggio Emilia, project:  {\it ``Discrete Methods in Combinatorial Geometry and Geometric Topology"}.

\end{document}